\documentclass[12pt,thmsa]{article}
\usepackage{sw20bams}

\input tcilatex
\begin{document}

\author{Steven Finch and Pascal Sebah}
\title{Residue of a Mod $5$ Euler Product}
\date{December 18, 2009}
\maketitle

\begin{abstract}
Consider the product of $(1-p^{-s})^{-4}$ over all primes $p\equiv 1\limfunc{%
mod}5$. We evaluate its residue at $s=1$ and compare with the corresponding
Mertens constant of Languasco \&\ Zaccagnini. We also count primitive
quintic Dirichlet characters mod $n$ and determine their average number as $%
n\rightarrow \infty $.
\end{abstract}

\footnotetext{
Copyright \copyright\ 2009 by Steven R. Finch. All rights reserved.}Let $%
\varphi (\ell )$ denote the number of positive integers $k\leq \ell $
satisfying $\gcd (k,\ell )=1$ and, for each such pair, 
\[
\rho (k,\ell )=\lim_{s\rightarrow 1}(s-1)\dprod\limits_{p\equiv k\limfunc{mod%
}\ell }\left( 1-\frac 1{p^s}\right) ^{-\varphi (\ell )} 
\]
where the product is over primes $p$. Clearly $\rho (1,2)=1/2$ since the
Euler product expression 
\[
\begin{array}{ccc}
\zeta (s)=\dprod\limits_p\left( 1-\dfrac 1{p^s}\right) ^{-1}, &  & \limfunc{%
Re}(s)>1
\end{array}
\]
for the Riemann zeta function has residue at $s=1$ equal to $1$. In the
event $\varphi (\ell )=2$, that is, $\ell \in \{3,4,6\}$, it is
straightforward to show that (section \ref{Euler346}) 
\[
\rho (1,3)=\frac{\sqrt{3}}{2\pi }\dprod\limits_{p\equiv 1\limfunc{mod}%
3}\left( 1-\frac 1{p^2}\right) ^{-1}, 
\]
\[
\rho (1,4)=\frac 1\pi \dprod\limits_{p\equiv 1\limfunc{mod}4}\left( 1-\frac
1{p^2}\right) ^{-1}, 
\]
\[
\rho (1,6)=\frac{\sqrt{3}}{2\pi }\dprod\limits_{p\equiv 1\limfunc{mod}%
6}\left( 1-\frac 1{p^2}\right) ^{-1}. 
\]
There are two outcomes corresponding to the event $\varphi (\ell )=2$. If $%
\ell \in \{8,12\}$, then (section \ref{Euler812}) 
\[
\rho (1,8)=\frac{2\ln \left( 1+\sqrt{2}\right) }{\pi ^2}\dprod\limits_{p%
\equiv 1\limfunc{mod}8}\left( 1-\frac 1{p^2}\right) ^{-2}, 
\]
\[
\rho (1,12)=\frac{3\ln \left( 2+\sqrt{3}\right) }{2\pi ^2}%
\dprod\limits_{p\equiv 1\limfunc{mod}12}\left( 1-\frac 1{p^2}\right) ^{-2}. 
\]
If $\ell \in \{5,10\}$, then (section \ref{Euler510}) 
\[
\rho (1,5)=\frac{5\ln \left( 2+\sqrt{5}\right) }{3\pi ^2}\dprod\limits_{p%
\equiv 1\limfunc{mod}5}\left( 1-\frac 1{p^2}\right) ^{-2}, 
\]
\[
\rho (1,10)=\frac{5\ln \left( 2+\sqrt{5}\right) }{3\pi ^2}%
\dprod\limits_{p\equiv 1\limfunc{mod}10}\left( 1-\frac 1{p^2}\right) ^{-2}. 
\]

We point out the equivalence of our computations with Languasco \&\
Zaccagnini \cite{LZ1, LZ2} via 
\[
e^{\gamma /\varphi (\ell )}\lim_{P\rightarrow \infty }\ln (P)^{1/\varphi
(\ell )}\dprod\limits\Sb p<P,  \\ p\equiv k\limfunc{mod}\ell  \endSb \left(
1-\frac 1p\right) =\lim_{s\rightarrow 1}\zeta (s)^{1/\varphi (\ell
)}\dprod\limits_{p\equiv k\limfunc{mod}\ell }\left( 1-\frac 1{p^s}\right) 
\]
for positive integers $k\leq \ell $ satisfying $\gcd (k,\ell )=1$, where $%
\gamma $ is the Euler-Mascheroni constant. Define the $(k,\ell )^{\text{th}}$
Mertens constant to be 
\[
\mu (k,\ell )=\lim_{P\rightarrow \infty }\ln (P)^{1/\varphi (\ell
)}\dprod\limits\Sb p<P,  \\ p\equiv k\limfunc{mod}\ell  \endSb \left(
1-\frac 1p\right) , 
\]
then the formula $\mu (1,1)=e^{-\gamma }$ is classical. From $e^{\gamma
/\varphi (\ell )}\mu (k,\ell )=\rho (k,\ell )^{-1/\varphi (\ell )}$, we have 
\[
\rho (k,\ell )=e^{-\gamma }\mu (k,\ell )^{-\varphi (\ell )}. 
\]
The cases $\ell =4,5,6,8,15,24$ appear explicitly in \cite{LZ1}; of the
first four of these, only their formula for $\ell =5$: 
\[
\rho (1,5)=\frac{50\ln \left( 2+\sqrt{5}\right) }{13\pi ^2\sqrt{5}}%
\dprod\limits_{p\equiv 1\limfunc{mod}5}\left( 1-\frac 1{p^4}\right)
^{-1}\dprod\limits_{p\equiv 4\limfunc{mod}5}\left( \frac{1+p^{-2}}{1-p^{-2}}%
\right) ^{-1} 
\]
looks symbolically different from our result. Together, the two formulas
yield 
\[
\dprod\limits_{p\equiv 1,4\limfunc{mod}5}\frac{p^2+1}{p^2-1}=\frac{6\sqrt{5}%
}{13} 
\]
which is known to be true. High-precision numerical values of $\mu (k,\ell )$
(what Languasco \&\ Zaccagnini call $C(\ell ,k)$) appear in \cite{LZ3}.

\section{Primes in Arithmetic Progression}

In the following, symbols $A$, $B$, $C$, $\ldots $ serve as placeholders for
nonzero values of a Dirichlet character. Superscripts on $A=A^1$, $B=B^1$, $%
C=C^1$, $\ldots $ suggest not exponentiation, but rather indicate the
associated roots of unity.

\subsection{\label{Euler346}Cases $\ell =3,4,6$}

Let $\chi _1$ denote the principal character mod $\ell $ and $\chi _2$
denote the nonprincipal character mod $\ell $. The corresponding L-series $%
L_1$ and $L_2$ satisfy 
\[
L_1(s)=\dprod\limits_p\left( 1-\frac{\chi _1(p)}{p^s}\right)
^{-1}=\dprod\limits\Sb p\equiv 1  \\ \limfunc{mod}\ell  \endSb \left(
1-\frac 1{p^s}\right) ^{-1}\cdot \dprod\limits\Sb p\equiv -1  \\ \limfunc{mod%
}\ell  \endSb \left( 1-\frac 1{p^s}\right) ^{-1}=AB, 
\]
\[
L_2(s)=\dprod\limits_p\left( 1-\frac{\chi _2(p)}{p^s}\right)
^{-1}=\dprod\limits\Sb p\equiv 1  \\ \limfunc{mod}\ell  \endSb \left(
1-\frac 1{p^s}\right) ^{-1}\cdot \dprod\limits\Sb p\equiv -1  \\ \limfunc{mod%
}\ell  \endSb \left( 1+\frac 1{p^s}\right) ^{-1}=AB^{-1}. 
\]
We have 
\[
L_1(2s)=\dprod\limits\Sb p\equiv 1  \\ \limfunc{mod}\ell  \endSb \left(
1-\frac 1{p^s}\right) ^{-1}\cdot \dprod\limits\Sb p\equiv 1  \\ \limfunc{mod}%
\ell  \endSb \left( 1+\frac 1{p^s}\right) ^{-1}\cdot \dprod\limits\Sb %
p\equiv -1  \\ \limfunc{mod}\ell  \endSb \left( 1-\frac 1{p^s}\right)
^{-1}\cdot \dprod\limits\Sb p\equiv -1  \\ \limfunc{mod}\ell  \endSb \left(
1+\frac 1{p^s}\right) ^{-1}=AA^{-1}BB^{-1}, 
\]
hence 
\[
g(s)=\frac{L_1(s)L_2(s)}{L_1(2s)}=\frac{AB\cdot AB^{-1}}{AA^{-1}BB^{-1}}%
=\frac A{A^{-1}}=\dprod\limits\Sb p\equiv 1  \\ \limfunc{mod}\ell  \endSb 
\frac{1+p^{-s}}{1-p^{-s}}. 
\]
From 
\[
\frac{1+x}{1-x}=(1-x)^{-2}(1-x^2), 
\]
we deduce that 
\[
g(s)=\dprod\limits\Sb p\equiv 1  \\ \limfunc{mod}\ell  \endSb \left( 1-\frac
1{p^s}\right) ^{-2}\left( 1-\frac 1{p^{2s}}\right) 
\]
and thus 
\[
\dprod\limits\Sb p\equiv 1  \\ \limfunc{mod}\ell  \endSb \left( 1-\frac
1{p^s}\right) ^{-2}=g(s)\dprod\limits\Sb p\equiv 1  \\ \limfunc{mod}\ell 
\endSb \left( 1-\frac 1{p^{2s}}\right) ^{-1}. 
\]
We also have 
\[
L_1(s)=\dprod\limits_{\gcd (p,\ell )=1}\left( 1-\frac 1{p^s}\right)
^{-1}=\dprod\limits_{p|\ell }\left( 1-\frac 1{p^s}\right) \zeta (s). 
\]
Therefore 
\begin{eqnarray*}
\lim_{s\rightarrow 1}\dprod\limits\Sb p\equiv 1  \\ \limfunc{mod}\ell 
\endSb \left( 1-\frac 1{p^s}\right) ^{-2}(s-1) &=&\left( \lim_{s\rightarrow
1}g(s)(s-1)\right) \dprod\limits\Sb p\equiv 1  \\ \limfunc{mod}\ell  \endSb %
\left( 1-\frac 1{p^2}\right) ^{-1} \\
\ &=&\left( \lim_{s\rightarrow 1}L_1(s)(s-1)\right) \frac{L_2(1)}{L_1(2)}%
\dprod\limits\Sb p\equiv 1  \\ \limfunc{mod}\ell  \endSb \left( 1-\frac
1{p^2}\right) ^{-1} \\
\ &=&\left( \lim_{s\rightarrow 1}\zeta (s)(s-1)\right) \dprod\limits_{p|\ell
}\frac{1-p^{-1}}{1-p^{-2}}\frac{L_2(1)}{\zeta (2)}\dprod\limits\Sb p\equiv 1 
\\ \limfunc{mod}\ell  \endSb \left( 1-\frac 1{p^2}\right) ^{-1} \\
\ &=&\frac 6{\pi ^2}\dprod\limits_{p|\ell }\frac p{p+1}L_2(1)\dprod\limits 
\Sb p\equiv 1  \\ \limfunc{mod}\ell  \endSb \left( 1-\frac 1{p^2}\right)
^{-1}.
\end{eqnarray*}
For example, if $\ell =4$, then $L_2(1)=\pi /4$ and the result simplifies: 
\[
\rho (1,4)=\frac 6{\pi ^2}\frac 23\frac \pi 4\dprod\limits\Sb p\equiv 1  \\ 
\limfunc{mod}4  \endSb \left( 1-\frac 1{p^2}\right) ^{-1}=\frac 1\pi
\dprod\limits\Sb p\equiv 1  \\ \limfunc{mod}4  \endSb \left( 1-\frac
1{p^2}\right) ^{-1}=\frac \pi {16\kappa ^2} 
\]
where $\kappa $ is the Landau-Ramanujan constant \cite{Fn}. Likewise $%
L_2(1)=\pi \sqrt{3}/9$ if $\ell =3$ and $L_2(1)=\pi \sqrt{3}/6$ if $\ell =6$%
, hence 
\[
\rho (1,3)=\frac 6{\pi ^2}\frac 34\frac{\pi \sqrt{3}}9\dprod\limits\Sb %
p\equiv 1  \\ \limfunc{mod}3  \endSb \left( 1-\frac 1{p^2}\right) ^{-1}=%
\frac{\sqrt{3}}{2\pi }\dprod\limits\Sb p\equiv 1  \\ \limfunc{mod}3  \endSb %
\left( 1-\frac 1{p^2}\right) ^{-1}, 
\]
\[
\rho (1,6)=\frac 6{\pi ^2}\frac 23\frac 34\frac{\pi \sqrt{3}}6\dprod\limits 
\Sb p\equiv 1  \\ \limfunc{mod}6  \endSb \left( 1-\frac 1{p^2}\right) ^{-1}=%
\frac{\sqrt{3}}{2\pi }\dprod\limits\Sb p\equiv 1  \\ \limfunc{mod}6  \endSb %
\left( 1-\frac 1{p^2}\right) ^{-1}. 
\]

\subsection{\label{Euler812}Cases $\ell =8,12$}

We have 
\begin{eqnarray*}
g(s) &=&\frac{L_1(s)L_2(s)L_3(s)L_4(s)}{L_1(2s)^2}=\frac{ABCD\cdot
ABC^{-1}D^{-1}\cdot AB^{-1}CD^{-1}\cdot AB^{-1}C^{-1}D}{%
(AA^{-1}BB^{-1}CC^{-1}DD^{-1})^2} \\
\ &=&\left( \frac A{A^{-1}}\right) ^2=\dprod\limits\Sb p\equiv 1  \\ 
\limfunc{mod}\ell  \endSb \left( \frac{1+p^{-s}}{1-p^{-s}}\right) ^2.
\end{eqnarray*}
From 
\[
\left( \frac{1+x}{1-x}\right) ^2=(1-x)^{-4}(1-x^2)^2, 
\]
we deduce that 
\[
g(s)=\dprod\limits\Sb p\equiv 1  \\ \limfunc{mod}\ell  \endSb \left( 1-\frac
1{p^s}\right) ^{-4}\left( 1-\frac 1{p^{2s}}\right) ^2 
\]
and thus 
\[
\dprod\limits\Sb p\equiv 1  \\ \limfunc{mod}\ell  \endSb \left( 1-\frac
1{p^s}\right) ^{-4}=g(s)\dprod\limits\Sb p\equiv 1  \\ \limfunc{mod}\ell 
\endSb \left( 1-\frac 1{p^{2s}}\right) ^{-2}. 
\]
Therefore 
\begin{eqnarray*}
\lim_{s\rightarrow 1}\dprod\limits\Sb p\equiv 1  \\ \limfunc{mod}\ell 
\endSb \left( 1-\frac 1{p^s}\right) ^{-4}(s-1) &=&\left( \lim_{s\rightarrow
1}g(s)(s-1)\right) \dprod\limits\Sb p\equiv 1  \\ \limfunc{mod}\ell  \endSb %
\left( 1-\frac 1{p^2}\right) ^{-2} \\
\ &=&\left( \lim_{s\rightarrow 1}L_1(s)(s-1)\right) \frac{L_2(1)L_3(1)L_4(1)%
}{L_1(2)^2}\dprod\limits\Sb p\equiv 1  \\ \limfunc{mod}\ell  \endSb \left(
1-\frac 1{p^2}\right) ^{-2} \\
\ &=&\left( \lim_{s\rightarrow 1}\zeta (s)(s-1)\right) \dprod\limits_{p|\ell
}\frac{1-p^{-1}}{(1-p^{-2})^2}\frac{L_2(1)L_3(1)L_4(1)}{\zeta (2)^2}%
\dprod\limits\Sb p\equiv 1  \\ \limfunc{mod}\ell  \endSb \left( 1-\frac
1{p^2}\right) ^{-2} \\
\ &=&\frac{36}{\pi ^4}\dprod\limits_{p|\ell }\frac{p^3}{(p+1)(p^2-1)}%
L_2(1)L_3(1)L_4(1)\dprod\limits\Sb p\equiv 1  \\ \limfunc{mod}\ell  \endSb %
\left( 1-\frac 1{p^2}\right) ^{-2}.
\end{eqnarray*}
If $\ell =8$, then 
\[
\begin{array}{ccccc}
L_2(1)=\dfrac{\pi \sqrt{2}}4, &  & L_3(1)=\dfrac \pi 4, &  & L_4(1)=\dfrac{%
\ln \left( 1+\sqrt{2}\right) }{\sqrt{2}}
\end{array}
\]
hence 
\[
\rho (1,8)=\frac{36}{\pi ^4}\frac 89\dfrac{\pi \sqrt{2}}4\dfrac \pi 4\dfrac{%
\ln \left( 1+\sqrt{2}\right) }{\sqrt{2}}\dprod\limits\Sb p\equiv 1  \\ 
\limfunc{mod}8  \endSb \left( 1-\frac 1{p^2}\right) ^{-2}=\frac{2\ln \left(
1+\sqrt{2}\right) }{\pi ^2}\dprod\limits\Sb p\equiv 1  \\ \limfunc{mod}8 
\endSb \left( 1-\frac 1{p^2}\right) ^{-2}. 
\]
If $\ell =12$, then 
\[
\begin{array}{ccccc}
L_2(1)=\dfrac \pi 3, &  & L_3(1)=\dfrac{\pi \sqrt{3}}6, &  & L_4(1)=\dfrac{%
\ln \left( 2+\sqrt{3}\right) }{\sqrt{3}}
\end{array}
\]
hence 
\[
\rho (1,12)=\frac{36}{\pi ^4}\frac 89\frac{27}{32}\dfrac \pi 3\dfrac{\pi 
\sqrt{3}}6\dfrac{\ln \left( 2+\sqrt{3}\right) }{\sqrt{3}}\dprod\limits\Sb %
p\equiv 1  \\ \limfunc{mod}12  \endSb \left( 1-\frac 1{p^2}\right) ^{-2}=%
\frac{3\ln \left( 2+\sqrt{3}\right) }{2\pi ^2}\dprod\limits\Sb p\equiv 1  \\ 
\limfunc{mod}12  \endSb \left( 1-\frac 1{p^2}\right) ^{-2}. 
\]

\subsection{\label{Euler510}Cases $\ell =5,10$}

Let $i$ denote the quartic root of unity. We have 
\begin{eqnarray*}
g(s) &=&\frac{L_1(s)L_2(s)L_3(s)L_4(s)}{L_1(2s)L_3(2s)}=\frac{ABCD\cdot
AB^iC^{-i}D^{-1}\cdot AB^{-1}C^{-1}D\cdot AB^{-i}C^iD^{-1}}{%
AA^{-1}BB^{-1}CC^{-1}DD^{-1}\cdot AA^{-1}B^iB^{-i}C^iC^{-i}DD^{-1}} \\
\ &=&\left( \frac A{A^{-1}}\right) ^2=\dprod\limits\Sb p\equiv 1  \\ 
\limfunc{mod}\ell  \endSb \left( \frac{1+p^{-s}}{1-p^{-s}}\right) ^2
\end{eqnarray*}
because 
\begin{eqnarray*}
L_3(2s) &=&\dprod\limits\Sb p\equiv 1  \\ \limfunc{mod}\ell  \endSb \left(
1-\frac 1{p^s}\right) ^{-1}\cdot \dprod\limits\Sb p\equiv 1  \\ \limfunc{mod}%
\ell  \endSb \left( 1+\frac 1{p^s}\right) ^{-1}\cdot \dprod\limits\Sb %
p\equiv 1+\ell /5  \\ \limfunc{mod}\ell  \endSb \left( 1-\frac i{p^s}\right)
^{-1}\cdot \dprod\limits\Sb p\equiv 1+\ell /5  \\ \limfunc{mod}\ell  \endSb %
\left( 1+\frac i{p^s}\right) ^{-1}\cdot \\
&&\;\;\;\cdot \dprod\limits\Sb p\equiv -1-\ell /5  \\ \limfunc{mod}\ell 
\endSb \left( 1-\frac i{p^s}\right) ^{-1}\cdot \dprod\limits\Sb p\equiv
-1-\ell /5  \\ \limfunc{mod}\ell  \endSb \left( 1+\frac i{p^s}\right)
^{-1}\cdot \dprod\limits\Sb p\equiv -1  \\ \limfunc{mod}\ell  \endSb \left(
1-\frac 1{p^s}\right) ^{-1}\cdot \dprod\limits\Sb p\equiv -1  \\ \limfunc{mod%
}\ell  \endSb \left( 1+\frac 1{p^s}\right) ^{-1} \\
&=&AA^{-1}B^iB^{-i}C^iC^{-i}DD^{-1}.
\end{eqnarray*}
By the same method as previously, 
\begin{eqnarray*}
\lim_{s\rightarrow 1}\dprod\limits\Sb p\equiv 1  \\ \limfunc{mod}\ell 
\endSb \left( 1-\frac 1{p^s}\right) ^{-4}(s-1) &=&\left( \lim_{s\rightarrow
1}L_1(s)(s-1)\right) \frac{L_2(1)L_3(1)L_4(1)}{L_1(2)L_3(2)}\dprod\limits\Sb %
p\equiv 1  \\ \limfunc{mod}\ell  \endSb \left( 1-\frac 1{p^2}\right) ^{-2} \\
\ &=&\left( \lim_{s\rightarrow 1}\zeta (s)(s-1)\right) \dprod\limits_{p|\ell
}\frac{1-p^{-1}}{1-p^{-2}}\frac{L_2(1)L_3(1)L_4(1)}{\zeta (2)L_3(2)}%
\dprod\limits\Sb p\equiv 1  \\ \limfunc{mod}\ell  \endSb \left( 1-\frac
1{p^2}\right) ^{-2} \\
\ &=&\frac 6{\pi ^2}\dprod\limits_{p|\ell }\frac p{p+1}\frac{%
L_2(1)L_3(1)L_4(1)}{L_3(2)}\dprod\limits\Sb p\equiv 1  \\ \limfunc{mod}\ell 
\endSb \left( 1-\frac 1{p^2}\right) ^{-2}.
\end{eqnarray*}
If $\ell =5$, then 
\[
\begin{array}{ccccc}
L_2(1)=2^{1/2}5^{-5/4}(3+4i)^{1/4}\pi , &  & L_3(1)=\dfrac{2\ln \left( 2+%
\sqrt{5}\right) }{3\sqrt{5}}, &  & L_4(1)=2^{1/2}5^{-5/4}(3-4i)^{1/4}\pi
\end{array}
\]
and $L_3(2)=4\pi ^2\sqrt{5}/125$, hence $L_2(1)L_4(1)/L_3(2)=\sqrt{5}/2$ and 
\[
\rho (1,5)=\frac 6{\pi ^2}\frac 56\frac{\sqrt{5}}2\dfrac{2\ln \left( 2+\sqrt{%
5}\right) }{3\sqrt{5}}\dprod\limits\Sb p\equiv 1  \\ \limfunc{mod}5  \endSb %
\left( 1-\frac 1{p^2}\right) ^{-2}=\frac{5\ln \left( 2+\sqrt{5}\right) }{%
3\pi ^2}\dprod\limits\Sb p\equiv 1  \\ \limfunc{mod}5  \endSb \left( 1-\frac
1{p^2}\right) ^{-2}. 
\]
If $\ell =10$, then 
\[
\begin{array}{ccccc}
L_2(1)=2^{-1/2}5^{-3/4}(3+4i)^{1/4}\pi , &  & L_3(1)=\dfrac{\ln \left( 2+%
\sqrt{5}\right) }{\sqrt{5}}, &  & L_4(1)=2^{-1/2}5^{-3/4}(3-4i)^{1/4}\pi
\end{array}
\]
and $L_3(2)=\pi ^2\sqrt{5}/25$, hence $L_2(1)L_4(1)/L_3(2)=\sqrt{5}/2$ and 
\[
\rho (1,10)=\frac 6{\pi ^2}\frac 23\frac 56\frac{\sqrt{5}}2\dfrac{\ln \left(
2+\sqrt{5}\right) }{\sqrt{5}}\dprod\limits\Sb p\equiv 1  \\ \limfunc{mod}10 
\endSb \left( 1-\frac 1{p^2}\right) ^{-2}=\frac{5\ln \left( 2+\sqrt{5}%
\right) }{3\pi ^2}\dprod\limits\Sb p\equiv 1  \\ \limfunc{mod}10  \endSb %
\left( 1-\frac 1{p^2}\right) ^{-2}. 
\]

\subsection{\label{Euler918}Cases $\ell =9,18$}

We will barely start the analysis, in the hope that someone else will finish
what we've begun. Let $\omega =\exp (i\pi /3)=(1+i\sqrt{3})/2$, the sextic
root of unity. We have 
\begin{eqnarray*}
h(s) &=&L_1(s)L_2(s)L_3(s)L_4(s)L_5(s)L_6(s) \\
\ &=&ABCDEF\cdot AB^\omega C^{\omega ^2}D^{-\omega ^2}E^{-\omega
}F^{-1}\cdot AB^{\omega ^2}C^{-\omega }D^{-\omega }E^{\omega ^2}F\cdot \\
&&\cdot AB^{-1}CD^{-1}EF^{-1}\cdot AB^{-\omega }C^{\omega ^2}D^{\omega
^2}E^{-\omega }F\cdot AB^{-\omega ^2}C^{-\omega }D^\omega E^{\omega ^2}F
\end{eqnarray*}
and must define an appropriate $g(s)=h(s)/\cdot $. L-values are more
complicated here; for example, when $\ell =9$, 
\[
L_3(1)=-\tfrac 23\omega ^{2/3}\left( -\omega \ln \left( \sin \left( \tfrac{%
2\pi }9\right) \right) +\omega ^2\ln \left( \cos \left( \tfrac \pi
{18}\right) \right) +\ln \left( \sin \left( \tfrac \pi 9\right) \right)
\right) . 
\]

\subsection{\label{Euler714}Cases $\ell =7,14$}

We have 
\begin{eqnarray*}
h(s) &=&L_1(s)L_2(s)L_3(s)L_4(s)L_5(s)L_6(s) \\
\ &=&ABCDEF\cdot AB^{\omega ^2}C^\omega D^{-\omega }E^{-\omega
^2}F^{-1}\cdot AB^{-\omega }C^{\omega ^2}D^{\omega ^2}E^{-\omega }F\cdot \\
&&\cdot ABC^{-1}DE^{-1}F^{-1}\cdot AB^{\omega ^2}C^{-\omega }D^{-\omega
}E^{\omega ^2}F\cdot AB^{-\omega }C^{-\omega ^2}D^{\omega ^2}E^\omega F^{-1}
\end{eqnarray*}
and must define an appropriate $g(s)=h(s)/\cdot $. L-values are more
complicated here; for example, when $\ell =7$, 
\[
L_5(1)=7^{-2/3}(-2-3\omega ^2)^{1/3}\left( -\omega \ln (y_1)+\omega ^2\ln
(y_2)+\ln (y_3)\right) 
\]
and $y_1<y_2<y_3$ are the (real)\ zeroes of $y^3-7y^2+14y-7$.

\section{Primitive Dirichlet Characters}

This is a follow-on to our earlier discussion \cite{FS, Fi} about asymptotic
enumeration of $\ell ^{\text{th}}$ order primitive Dirichlet characters mod $%
n$. Let $b_\ell (n)$ denote the count of such characters. There exists a
constant $0<K_\ell <\infty $ such that \cite{FMS} 
\[
\dsum\limits_{n\leq N}b_\ell (n)\sim K_\ell N\ln (N)^{d(\ell )-2} 
\]
as $N\rightarrow \infty $, where $d(\ell )$ is the number of divisors of $%
\ell $. Special cases of this result for $\ell =2,3,4,8$ were examined in 
\cite{Fi}. Leading coefficients are known to be $K_2=6/\pi ^2$, 
\[
K_3=\frac{11\sqrt{3}}{18\pi }\dprod\limits_{p\equiv 1\limfunc{mod}3}\left(
1-\frac 2{p(p+1)}\right) , 
\]
\[
K_4=\frac 7{\pi ^3}\dprod\limits_{p\equiv 1\limfunc{mod}4}\left( 1-\frac
1{p^2}\right) ^{-1}\left( 1-\frac{5p-3}{p^2(p+1)}\right) , 
\]
\begin{eqnarray*}
K_8 &=&\frac 1{2!}\frac{16\ln (1+\sqrt{2})}{\pi ^5}\dprod\limits_{p\equiv 5%
\limfunc{mod}8}\left( 1-\dfrac 1{p^2}\right) ^{-1}\left( 1-\dfrac{5p-3}{%
p^2(p+1)}\right) \cdot \\
&&\ \cdot \dprod\limits_{p\equiv 1\limfunc{mod}8}\left( 1-\dfrac
1{p^2}\right) ^{-3}\left( 1-\dfrac{27p^5-85p^4+125p^3-99p^2+41p-7}{p^6(p+1)}%
\right) .
\end{eqnarray*}
The case for $\ell =6$ served as an illustration of a general method in \cite
{FMS}: 
\[
K_6=\frac 1{2!}\frac 6{\pi ^4}\dprod\limits_{p\equiv 1\limfunc{mod}6}\left(
1-\frac 1{p^2}\right) ^{-2}\left( 1-\frac{14p^3-26p^2+19p-5}{p^4(p+1)}%
\right) 
\]
and we will review its derivation here. A ``fairly roundabout'' procedure
for calculating required Euler product residues in \cite{Fi} is now avoided.
We will also study the case for $\ell =5$; the cases for $\ell =7,9$ seem
more difficult.

Since 
\[
\left( 1-\frac 1{p^2}\right) ^{-1}\left( 1-\frac{5p-3}{p^2(p+1)}\right)
=1-\frac 4{(p+1)^2}, 
\]
\[
\left( 1-\frac 1{p^2}\right) ^{-2}\left( 1-\frac{14p^3-26p^2+19p-5}{p^4(p+1)}%
\right) =1-\frac{12p-4}{(p+1)^3}, 
\]
\[
\left( 1-\dfrac 1{p^2}\right) ^{-3}\left( 1-\dfrac{%
27p^5-85p^4+125p^3-99p^2+41p-7}{p^6(p+1)}\right) =1-\frac{24p^2-16p+8}{%
(p+1)^4} 
\]
alternative, more compact expressions for $K_4$, $K_6$, $K_8$ are available
as well \cite{FMS}.

\subsection{Sextic Characters}

Subscripts are omitted for simplicity. When $\ell =6$, we have 
\[
\begin{array}{ccc}
b(2^r)=\left\{ 
\begin{array}{lll}
0 &  & \text{if }r=1, \\ 
1 &  & \text{if }r=2, \\ 
2 &  & \text{if }r=3, \\ 
0 &  & \text{if }r\geq 4,
\end{array}
\right. &  & b(3^r)=\left\{ 
\begin{array}{lll}
1 &  & \text{if }r=1, \\ 
4 &  & \text{if }r=2, \\ 
0 &  & \text{if }r\geq 3,
\end{array}
\right.
\end{array}
\]
\[
b(p^r)=\left\{ 
\begin{array}{lll}
1 &  & \text{if }r=1\text{ \& }p\equiv 5\limfunc{mod}6, \\ 
5 &  & \text{if }r=1\text{ \& }p\equiv 1\limfunc{mod}6 \\ 
0 &  & \text{otherwise}
\end{array}
\right. 
\]
for prime $p\geq 5$ and $r\geq 1$, hence 
\begin{eqnarray*}
\dsum\limits_{n=1}^\infty \frac{b(n)}{n^s} &=&\left( 1+\frac 1{2^{2s}}+\frac
2{2^{3s}}\right) \left( 1+\frac 1{3^s}+\frac 4{3^{2s}}\right) \dprod\limits 
\Sb p\equiv 5  \\ \limfunc{mod}6  \endSb \left( 1+\frac 1{p^s}\right) \cdot
\dprod\limits\Sb p\equiv 1  \\ \limfunc{mod}6  \endSb \left( 1+\frac
5{p^s}\right) \\
\ &=&\left( 1+\frac 1{2^{2s}}+\frac 2{2^{3s}}\right) \left( 1+\frac
1{3^s}+\frac 4{3^{2s}}\right) \dprod\limits\Sb p\equiv 5  \\ \limfunc{mod}6 
\endSb \left( 1-\frac 1{p^s}\right) ^{-1}\left( 1-\frac 1{p^{2s}}\right)
\cdot \\
&&\ \cdot \dprod\limits\Sb p\equiv 1  \\ \limfunc{mod}6  \endSb \left(
1-\frac 1{p^s}\right) ^{-5}\left( 1-\frac 1{p^{2s}}\right) \left( 1-\frac{%
14p^{3s}-26p^{2s}+19p^s-5}{p^{4s}(p^s+1)}\right) \\
\ &=&\left( 1+\frac 1{2^{2s}}+\frac 2{2^{3s}}\right) \left( 1+\frac
1{3^s}+\frac 4{3^{2s}}\right) \left( 1-\frac 1{2^s}\right) \left( 1-\frac
1{3^s}\right) \zeta (s)\left( 1-\frac 1{2^{2s}}\right) ^{-1}\cdot \\
&&\ \cdot \left( 1-\frac 1{3^{2s}}\right) ^{-1}\zeta (2s)^{-1}\dprod\limits 
\Sb p\equiv 1  \\ \limfunc{mod}6  \endSb \left( 1-\frac 1{p^s}\right)
^{-4}\left( 1-\frac{14p^{3s}-26p^{2s}+19p^s-5}{p^{4s}(p^s+1)}\right)
\end{eqnarray*}
and (from section \ref{Euler346}) 
\[
\dprod\limits_{p\equiv 1\limfunc{mod}6}\left( 1-\frac 1{p^s}\right)
^{-2}\sim \zeta (s)\frac{\sqrt{3}}{2\pi }\dprod\limits_{p\equiv 1\limfunc{mod%
}6}\left( 1-\frac 1{p^2}\right) ^{-1} 
\]
as $s\rightarrow 1$. The overall exponent of $\zeta (s)\,$ is $3=d(6)-1$,
consistent with \cite{FMS}; also the expression for $K_6$ follows by the
Selberg-Delange method \cite{FS}.

\subsection{Quintic Characters}

When $\ell =5$, we have 
\[
\begin{array}{ccc}
b(5^r)=\left\{ 
\begin{array}{lll}
4 &  & \text{if }r=2, \\ 
0 &  & \text{otherwise},
\end{array}
\right. &  & b(p^r)=\left\{ 
\begin{array}{lll}
4 &  & \text{if }r=1\text{ \& }p\equiv 1\limfunc{mod}5, \\ 
0 &  & \text{otherwise}
\end{array}
\right.
\end{array}
\]
for prime $p\neq 5$ and $r\geq 1$, hence 
\begin{eqnarray*}
\dsum\limits_{n=1}^\infty \frac{b(n)}{n^s} &=&\left( 1+\frac
4{5^{2s}}\right) \dprod\limits\Sb p\equiv 1  \\ \limfunc{mod}5  \endSb %
\left( 1+\frac 4{p^s}\right) \\
\ &=&\left( 1+\frac 4{5^{2s}}\right) \dprod\limits\Sb p\equiv 1  \\ \limfunc{%
mod}5  \endSb \left( 1-\frac 1{p^s}\right) ^{-4}\left( 1-\frac
1{p^{2s}}\right) \left( 1-\frac{9p^{2s}-11p^s+4}{p^{3s}(p^s+1)}\right)
\end{eqnarray*}
and (from section \ref{Euler510}) 
\[
\dprod\limits_{p\equiv 1\limfunc{mod}5}\left( 1-\frac 1{p^s}\right)
^{-4}\sim \zeta (s)\frac{5\ln \left( 2+\sqrt{5}\right) }{3\pi ^2}%
\dprod\limits_{p\equiv 1\limfunc{mod}5}\left( 1-\frac 1{p^2}\right) ^{-2} 
\]
as $s\rightarrow 1$. The overall exponent of $\zeta (s)\,$ is $1=d(5)-1$,
consistent with \cite{FMS}; also 
\begin{eqnarray*}
K_5 &=&\frac{29\ln \left( 2+\sqrt{5}\right) }{15\pi ^2}\dprod\limits\Sb %
p\equiv 1  \\ \limfunc{mod}5  \endSb \left( 1-\frac 1{p^2}\right)
^{-1}\left( 1-\frac{9p^2-11p+4}{p^3(p+1)}\right) \\
&=&\frac{29\ln \left( 2+\sqrt{5}\right) }{15\pi ^2}\dprod\limits\Sb p\equiv
1  \\ \limfunc{mod}5  \endSb \left( 1-\frac{8p-4}{p(p+1)^2}\right)
\end{eqnarray*}
by the Selberg-Delange method \cite{FS}. This is the average number of
primitive quintic Dirichlet characters mod $n$ as $n\rightarrow \infty $.

Starting points for numerical verification of the value $K_5$ might begin
with 
\[
\mu (1,5)=1.22523843853908458005760977474922052754059550939164... 
\]
from \cite{LZ3} and 
\[
K_5=\frac{29}{25}\dprod\limits_p\left( 1+\frac{\gcd (5,p-1)}{p-1}\right)
\left( 1-\frac 1p\right) ^2 
\]
as a specialization of (slowly convergent) general formulas in \cite{FMS}.

\subsection{Septic Characters}

When $\ell =7$, we have 
\[
\begin{array}{ccc}
b(7^r)=\left\{ 
\begin{array}{lll}
6 &  & \text{if }r=2, \\ 
0 &  & \text{otherwise},
\end{array}
\right. &  & b(p^r)=\left\{ 
\begin{array}{lll}
6 &  & \text{if }r=1\text{ \& }p\equiv 1\limfunc{mod}7, \\ 
0 &  & \text{otherwise}
\end{array}
\right.
\end{array}
\]
for prime $p\neq 7$ and $r\geq 1$, hence 
\begin{eqnarray*}
\dsum\limits_{n=1}^\infty \frac{b(n)}{n^s} &=&\left( 1+\frac
6{7^{2s}}\right) \dprod\limits\Sb p\equiv 1  \\ \limfunc{mod}7  \endSb %
\left( 1+\frac 6{p^s}\right) \\
\ &=&\left( 1+\frac 6{7^{2s}}\right) \dprod\limits\Sb p\equiv 1  \\ \limfunc{%
mod}7  \endSb \left( 1-\frac 1{p^s}\right) ^{-6}\left( 1-\frac
1{p^{2s}}\right) \left( 1-\frac{20p^{4s}-50p^{3s}+55p^{2s}-29p^s+6}{%
p^{5s}(p^s+1)}\right) .
\end{eqnarray*}
An infinite product formulation for $K_7$, akin to the others, awaits the
resolution of $\rho (1,7)$ (section \ref{Euler714}).

Numerical calculation of the value $K_7$ might begin with 
\[
\mu (1,7)=1.20435271605501440413126997959392601183676589049086... 
\]
from \cite{LZ3} and 
\[
K_7=\frac{55}{49}\dprod\limits_p\left( 1+\frac{\gcd (7,p-1)}{p-1}\right)
\left( 1-\frac 1p\right) ^2 
\]
from \cite{FMS}.

\subsection{Nonic Characters}

When $\ell =9$, we have 
\[
\begin{array}{ccc}
b(3^r)=\left\{ 
\begin{array}{lll}
2 &  & \text{if }r=2, \\ 
6 &  & \text{if }r=3, \\ 
0 &  & \text{if }r\geq 4,
\end{array}
\right. &  & b(p^r)=\left\{ 
\begin{array}{lll}
2 &  & \text{if }r=1\text{ \& }p\equiv 4,7\limfunc{mod}9, \\ 
8 &  & \text{if }r=1\text{ \& }p\equiv 1\limfunc{mod}9 \\ 
0 &  & \text{otherwise}
\end{array}
\right.
\end{array}
\]
for prime $p\neq 3$ and $r\geq 1$, hence 
\begin{eqnarray*}
\dsum\limits_{n=1}^\infty \frac{b(n)}{n^s} &=&\left( 1+\frac 2{3^{2s}}+\frac
6{3^{3s}}\right) \dprod\limits\Sb p\equiv 4,7  \\ \limfunc{mod}9  \endSb %
\left( 1+\frac 2{p^s}\right) \cdot \dprod\limits\Sb p\equiv 1  \\ \limfunc{%
mod}9  \endSb \left( 1+\frac 8{p^s}\right) \\
\ &=&\left( 1+\frac 2{3^{2s}}+\frac 6{3^{3s}}\right) \dprod\limits\Sb %
p\equiv 4,7  \\ \limfunc{mod}9  \endSb \left( 1-\frac 1{p^s}\right)
^{-2}\left( 1-\frac 1{p^{2s}}\right) \left( 1-\frac 2{p^s(p^s+1)}\right)
\cdot \\
&&\cdot \dprod\limits\Sb p\equiv 1  \\ \limfunc{mod}9  \endSb \left( 1-\frac
1{p^s}\right) ^{-8}\left( 1-\frac 1{p^{2s}}\right) \left( 1-\frac{%
35p^{6s}-133p^{5s}+245p^{4s}-259p^{3s}+161p^{2s}-55p^s+8}{p^{7s}(p+1)}%
\right) .
\end{eqnarray*}
An infinite product formulation for $K_9$, akin to the others, awaits the
resolution of not just $\rho (1,9)$ (section \ref{Euler918}), but also $\rho
(4,9)$ and $\rho (7,9)$.

Numerical calculation of the value $K_9$ might begin with 
\[
\mu (1,9)=1.17384958686544919027013946839197396049956269282192..., 
\]
\[
\mu (4,9)=1.13360386133436932499173359590759623742339637224206..., 
\]
\[
\mu (7,9)=1.05470661565485874510828199884014910243407287242835... 
\]
from \cite{LZ3} and 
\[
K_9=\frac{13}9\dprod\limits_p\left( 1+\frac{\gcd (9,p-1)}{p-1}\right) \left(
1-\frac 1p\right) ^3 
\]
from \cite{FMS}.

\section{Acknowledgement}

Portions of this paper were written before we met Greg Martin -- the general
formulas and theoretical rigor in \cite{FMS} are largely due to him -- thus
we wish to publicly express our appreciation to him for a very enjoyable and
rewarding collaboration.

\end{document}